\documentclass[11pt]{article}
                                                                                        
\newcommand{\comment}[1]{}

\newcommand{\reals}{\mbox{$\mathbb R$}}

\newcommand{\nats}{\mbox{$\mathbb N$}}

\newcommand{\power}{\mbox{$\mathbb P$}}
\newcommand{\floor}[1]{\left\lfloor #1 \right\rfloor}
\newcommand{\ceil}[1]{\left\lceil #1 \right\rceil}

\newcommand{\op}{{\mathop{\mathrm{op}}\nolimits}}
                                                                                        
\usepackage{amssymb}
                                                                                        
\def\squarebox#1{\hbox to #1{\hfill\vbox to #1{\vfill}}}
\def\qed{\hspace*{\fill}
        \vbox{\hrule\hbox{\vrule\squarebox{.667em}\vrule}\hrule}\smallskip}
\newenvironment{proof}{\begin{trivlist}
  \item[\hspace{\labelsep}{\em\noindent Proof.~}]
  }{\qed\end{trivlist}}
                                                                                        
\newtheorem{lemma}{Lemma}[section]
\newtheorem{theorem}[lemma]{Theorem}
\newtheorem{corollary}[lemma]{Corollary}

\newtheorem{observation}[lemma]{Observation}
\newtheorem{definition}[lemma]{Definition}

\def\squareforqed{\hbox{\rlap{$\sqcap$}$\sqcup$}}
\def\qed{\ifmmode\squareforqed\else{\unskip\nobreak\hfil
\penalty50\hskip1em\null\nobreak\hfil\squareforqed
\parfillskip=0pt\finalhyphendemerits=0\endgraf}\fi}
                                                                                        
\newlength{\tablength}
\newlength{\spacelength}
\newcommand{\tabstar}{\hspace*{\tablength}}
\newcommand{\spacestar}{\hspace*{\spacelength}}
\def\obeytabs{\catcode`\^^I=\active}
{\obeytabs\global\let^^I=\tabstar}
{\obeyspaces\global\let =\spacestar}
\newenvironment{display}{\begingroup\obeylines\obeyspaces\obeytabs}{\endgroup}
\newenvironment{prog}{\begin{display}\parskip0pt\sf}{\end{display}}

\title{On multipartite posets}

\author{
{\sl Geir Agnarsson}\thanks{ Department of Mathematical Sciences;
George Mason University, MS 3F2; 4400 University Drive; 
Fairfax, VA -- 22030; {\tt geir@math.gmu.edu}} 
}

\date{}

\begin{document}

\maketitle

\begin{abstract}
A poset $\mathbf{P} = (X,\preceq)$ is {\em $m$-partite} if $X$ has a partition
$X = X_1 \cup \cdots \cup X_m$ such that (1) each $X_i$ forms an 
antichain in $\mathbf{P}$, and (2) $x\prec y$ implies $x\in X_i$ and
$y\in X_j$ where $i<j$. In this article we derive a tight asymptotic 
upper bound on the order dimension of $m$-partite posets in terms of $m$ 
and their bipartite sub-posets in a constructive and elementary way. 

\vspace{3 mm}



\vspace{2 mm}

\end{abstract}

\section{Introduction}
\label{sec:intro}

The purpose of this article is to derive an asymptotically tight
upper bound for the dimension of multipartite posets in terms of their
number of parts and their bipartite sub-posets. 
Precise definitions of terms will be given later in Section~\ref{sec:defs}.
This work was partly inspired by a question asked by 
Reinhard Laubenbacher~\cite{Rein-privat} which casually can be phrased as follows: 
``For a given collection of posets, form a new poset
by stacking them together, putting one on top of the other. Is
it possible to bound the order dimension of the newly
formed poset in terms of the order dimension of the given posets?''
Laubenbacher's motivation were posets that appeared in the following manner:
When finitely many agents $A_1,\ldots, A_n$ are investigated 
over discrete times $t = 0,1,\ldots,m$, one obtains a poset 
consisting of the $n(m+1)$ elements $A_i(t)$, where 
a directed edge from $A_i(t)$ down to $A_j(t+1)$ is present if, and only
if, agent $A_i$ has influenced agent $A_j$ during the time interval
from $t$ to $t+1$. This resulting induced poset is sometimes called 
the {\em influence poset} among the agents. Here we have $m+1$ parts
of the influence poset, one part $X_t = \{A_1(t),\ldots,A_n(t)\}$ 
for each time $t = 0,1,\ldots,m$.

Other more classical posets can also be viewed as stacked sub-posets,
one on top of the other: If ${\cal{F}}_P$ is the face lattice of an 
$n$-dimensional polytope $P$ and ${\cal{F}}_P(i,i+1)$ is the height-2
sub-poset
of ${\cal{F}}_P$ consisting of the $i$ and $(i+1)$-dimensional faces of $P$,
then ${\cal{F}}_P$ can be thought of being formed by stacking 
${\cal{F}}_P(i,i+1)$ on top of ${\cal{F}}_P(i-1,i)$ for each $i=0,1,\ldots,n$.
In this case the stacking appears naturally since ${\cal{F}}_P$ is
a {\em graded poset} provided with a grading function into 
the nonnegative integers, that maps each face of $P$ (i.e.~each element of 
the poset ${\cal{F}}_P$) to its dimension.
(For more on graded posets see~\cite{Stanley1} and~\cite{Schroeder}.)
Determining the order dimension of face lattices of convex polytopes
is hard. Some partial yet interesting results in this direction
appear in~\cite{Reuter} and later in~\cite{Bright-Trott}.
Of particular interest in the literature is the face lattice of the 
standard $n$-simplex when viewed as the subset lattice of $\{1,\ldots,n\}$.
If we let $[n] = \{1,\ldots,n\}$ and 
${[n]\choose k}$ denote all the $k$-element subsets of
$[n]$, then the power set $\power([n])$ of all subsets of 
$[n]$ can be partitioned into $n+1$ disjoint sets 
$\power([n]) = {[n]\choose 0}\cup {[n]\choose 1}\cup\cdots\cup{[n]\choose n}$.
For $0\leq k_1<k_2\leq n$ denote the poset on 
${[n]\choose k_1}\cup {[n]\choose k_2}$ induced by
inclusion by $\mathbf{P}(k_1,k_2;n)$. Hence, as a poset,
$\power([n])$ can 
be thought of being formed by stacking the $n$ 
posets $\mathbf{P}(i,i+1;n)$ for $i\in \{0,1,\ldots,n-1\}$
one on top of the other.
Investigating the order dimension $\dim(k_1,k_2;n)$ 
of such sub-posets $\mathbf{P}(k_1,k_2;n)$
of $\power([n])$ for $1\leq k_1<k_2\leq n-1$ is currently an active area
of research, in particular the investigation of $\dim(1,k;n)$,
the order dimension of the poset $\mathbf{P}(1,k;n)$.
We list and briefly discuss a few related and celebrated results of this 
ongoing investigation: 

In~\cite{Dushnik} an explicit formula for $\dim(1,k;n)$ is
given when $2\sqrt{n}-2\leq k < n-1$ and in \cite{Trotter} 
the exact values of $\dim(1,2;n)$ are given for $2\leq n\leq 13$. 
In~\cite{Hurlbert} and~\cite{Furedi}
the exact values for $\dim(2,n-2;n)$ and $\dim(k,n-k;n)$ are given, 
provided that certain conditions hold for $k$ and $n$. 
In~\cite{Spencer} the asymptotic behavior of $\dim(1,k;n)$ is
given as a function of $n$ when $k$ is considered fixed.
Finally, in~\cite{Serkan-Morris} a direct method to determine 
$\dim(1,2;n)$ for each $n$ is given. 
Hence, the case $k=2$ for determining $\dim(1,k;n)$ is the only 
case which can be considered completely solved.
In~\cite{Kierstead} however, it is shown by contradiction that 
$\dim(1,\log n,n)=\Omega(\log^3 n/\log(\log n))$.
In addition, all the upper bounds derived there are proved 
by explicit construction, which therefore is also
an effective method in providing bounds for order dimensions.

In what follows we will discuss a class of posets that will
include the class of graded posets and the posets obtained
by such ``stacking'' as mentioned above in an ad hoc manner.
Our methods will be constructive and combinatorially elementary. 
In Section~\ref{sec:defs} we introduce our notation, state
our definitions in a precise manner and dispatch some basic properties. 
In the last Section~\ref{sec:multi-main} we state and prove our main result
of this article.

\section{Definitions and basic properties}
\label{sec:defs}

By a {\em poset} $\mathbf{P}$ we will always mean an ordered
tuple $\mathbf{P} = (X,\preceq)$ where $\preceq$ is a 
reflexive, antisymmetric and transitive binary 
relation on $X$. Unless otherwise stated $X$ is 
always assumed to be a finite set. We will for the most
part try to be consistent with the standard notation from~\cite{Trottbok}.
In particular, if two elements $x,y\in X$ are incomparable in $\mathbf{P}$,
then we write $x\parallel y$. By $\min(\mathbf{P})$ and $\max(\mathbf{P})$
we mean the set of minimal and maximal elements of $\mathbf{P}$ 
respectively. As originally defined in~\cite{Dush-Mill} and as stated
in~\cite{Trottbok}, the {\em order dimension} of 
$\mathbf{P} = (X,\preceq)$, denoted by $\dim(\mathbf{P})$, 
is the least number $d\in\nats$ of linear extensions 
$\preceq_1,\ldots, \preceq_d$ of $\preceq$ that {\em realize} $\preceq$.
That is, for $x,y\in X$ we have $x\preceq y$ in $\mathbf{P}$ iff $x\preceq_i y$ for 
all $i\in [d]$. 

Recall that for $n\in \nats$, any collection $S$ of points in 
the $n$-dimensional Euclidean space
${\reals}^n$ naturally forms a poset $(S,\preceq_E)$ by
$\tilde{x}\preceq_E\tilde{y} \Leftrightarrow x_i\leq y_i
\mbox{ for each } i\in [n]$, for any $\tilde{x} = (x_1,\ldots,x_n)$ 
and $\tilde{y} = (y_1,\ldots,y_n)$
from $S$. With this in mind we have that the order dimension
$\dim(\mathbf{P})$ of a poset $\mathbf{P} = (X,\preceq)$ 
is the least $d\in\nats$ such that there 
is an injective homomorphism $\phi : \mathbf{P}\rightarrow {\reals}^d$ satisfying
$x\preceq y \Leftrightarrow \phi(x)\preceq_E \phi(y)$ for all
$x,y\in X$. Hence the words order {\em dimension}.
Determining the exact value of the order dimension of a poset
is a hard computational problem. Even when we restrict to height-2
posets, the problem of computing their order
dimensions is NP-complete~\cite{Yanna}.

Recall that the ``standard example'' $\mathbf{S}_{n}$
from \cite{Dush-Mill} and \cite[p.~12]{Trottbok} is
a poset $\mathbf{S}_{n} = (A\cup B,\preceq)$ where
$A = \{a_1,\ldots,a_n\}$ and $B = \{b_1,\ldots,b_n\}$ are disjoint and
$a_i\prec b_j$ if, and only if, $i\neq j$. Here $S_n$ is a height-2 
poset on $2n$ elements with order dimension of $n$. 
By adding an element $c_{i\/j}$ between $a_i$ and $b_j$
for each $i\neq j$, so $a_i\prec c_{i\/j}\prec b_j$, we obtain
a poset $\mathbf{P} = (A\cup C\cup B,\preceq)$ on $n(n+1)$ 
elements induced by the $n(n-1)$ relations $a_i\prec c_{i\/j}\prec b_j$.
Clearly, the sub-poset induced by $A\cup B$ is the standard example,
so $\dim(\mathbf{P})\geq n$. However, $\mathbf{P}$ is obtained
by stacking the sub-poset induced by $B\cup C$ on top of the one
induced by $A\cup C$, each of which has the order dimension 2.
From this we obtain the following trivial but noteworthy observation.
\begin{observation}
\label{obs:no-func}
There is no function $f : \nats\times\nats\rightarrow \nats$
such that $\dim(\mathbf{P}) \leq f(\dim(\mathbf{P}_1),\dim(\mathbf{P}_2))$
holds in general for all posets 
$\mathbf{P}$, which are induced by two sub-posets 
$\mathbf{P}_1$ and $\mathbf{P}_2$ with
$\min(\mathbf{P}_1) = \max(\mathbf{P}_2)$.
\end{observation}
Although this answers the initial motivating question of Laubenbacher
from Section~\ref{sec:intro} in the negative, it does prompt us to bound
the order dimension in terms of other sub-posets.

The following lemma is a direct consequence of
the interpolation property for posets
and the fact that each poset has a linear extension.
\begin{lemma}
\label{lmm:extension}
Let $\mathbf{P}$ be a poset and $\mathbf{P}'$ be an induced
sub-poset of $\mathbf{P}$. Then any linear extension $\mathbf{L}'$
of $\mathbf{P}'$ can be extended to a linear extension $\mathbf{L}$
of $\mathbf{P}$.
\end{lemma}
Recall that a {\em bipartite} poset is an ordered triple
$\mathbf{P} = (X,Y;\preceq)$ where $X$ and $Y$ are disjoint
and $x\prec y$ implies that $x\in X$ and $y\in Y$. This can be generalized.
\begin{definition}
\label{def:multipartite}
Let $m\geq 2$ be an integer and
$X_1,\ldots, X_m$ be disjoint nonempty sets. We call 
$\mathbf{P} = (X_1,\ldots,X_m;\preceq)$ an {\em $m$-partite poset}
if $\preceq$ is a partial order on $X = X_1\cup\cdots\cup X_m$
such that (1) each $X_i$ forms an antichain w.r.t.~$\preceq$,
and (2) $x\prec y$ implies $x\in X_i$ and $y\in X_j$ where
$i,j\in [m]$ and $i<j$.
If $\mathbf{P}$ is $m$-partite for some $m$, then $\mathbf{P}$ is
a {\em multipartite poset}.
\end{definition}
Clearly, each $m$-partite poset $\mathbf{P}$ yields
its underlying poset $\mathbf{P}^{\circ} = (X_1\cup\cdots\cup X_m,\preceq)$ 
by ignoring the partition. The order dimension of 
$\mathbf{P}$ is then defined to be that of $\mathbf{P}^{\circ}$.

\section{Multipartite posets}
\label{sec:multi-main}

Let $\mathbf{P} = (X_1,\ldots,X_m;\preceq)$
be an $m$-partite poset and $\mathbf{P}_{i,j}$
be the bipartite sub-poset of $\mathbf{P}$ 
induced by $X_i\cup X_j$ for each $i<j$ with $i,j\in [m]$.
By Observation~\ref{obs:no-func},
we cannot hope to express $\dim(\mathbf{P})$ in terms of the 
$\dim(\mathbf{P}_{i,i+1})$'s for $i\in[m-1]$, the order dimensions
of these consecutive layers in $\mathbf{P}$. More is needed.

For each $i,j\in [m]$ with $i<j$ let 
$d_{i,j} = \dim(\mathbf{P}_{i,j})$ 
and ${\mathcal{L}}_{i,j}$ be a 
collection of $d_{i,j}$ linear orders on $X_i\cup X_j$ 
realizing $\mathbf{P}_{i,j}$.
By Lemma~\ref{lmm:extension} there is a set 
${\mathcal{L}}^{*}_{i,j}$ of $d_{i,j}$ linear orders extending
$\mathbf{P}$ and each linear order in ${\mathcal{L}}_{i,j}$.
By considering both cases of $x\parallel y$, where
$x,y\in X_i$ for some $i$ on one hand, and 
$x\in X_i$, $y\in X_j$ for some $i\neq j$ on the other,
we can see that $\mathcal{R} = \bigcup_{i<j}{\mathcal{L}}^{*}_{i,j}$
realizes $\mathbf{P}$. This shows that we can 
bound $\dim(\mathbf{P})$ in terms of the 
$\dim(\mathbf{P}_{i,j})$'s. We summarize in the following.
\begin{observation}
\label{obs:multi-upper}
For a multipartite poset $\mathbf{P}=(X_1,\ldots,X_m;\preceq)$
we have
\[
\dim(\mathbf{P})\leq 
\sum_{i<j}\dim(\mathbf{P}_{i,j}).
\]
\end{observation}
For an $m$-partite poset $\mathbf{P}$ let 
$B(\mathbf{P}) = \max_{i<j}\{\dim(\mathbf{P}_{i,j})\}$.
Since there are ${m\choose 2} = m(m-1)/2$ posets 
$\mathbf{P}_{i,j}$ we obtain
\[
B(\mathbf{P}) \leq \dim(\mathbf{P})
\leq\frac{m(m-1)}{2}B(\mathbf{P}),
\]
and hence for a fixed $m$, we have 
$\dim(\mathbf{P}) = \Theta(B(\mathbf{P}))$.
This can be reduced by a factor of $1/2$ in the following theorem.
\begin{theorem}
\label{thm:better-bounds}
For a multipartite poset $\mathbf{P} = (X_1,\ldots,X_m;\preceq)$
we have  
\[
\dim(\mathbf{P})
\leq \floor{\frac{(m-1)(m+3)}{4}}B(\mathbf{P}).
\]
\end{theorem}
\begin{proof}
Note that if $i_1<j_1<i_2<j_2<\cdots <i_{\ell}<j_{\ell}$
are indices from $[m]$ and $\mathbf{L}_k$ 
is a linear extension of $\mathbf{P}_{i_k,j_k}$,
then a linear extension of $\mathbf{P}$ that includes  
$\mathbf{L}_1\prec \mathbf{L}_2\prec \cdots\prec \mathbf{L}_{\ell}$
extends $\mathbf{P}$ and each of the $\mathbf{L}_k$.
In this way we can find $2\cdot B(\mathbf{P})$
linear orders extending
$\mathbf{P}$ and each $\mathbf{L}_{i,j}\in {\mathcal{L}}_{i,j}$, 
where $i+1 = j$. In general, for each $k\leq \floor{(m+1)/2}$ there are 
$k\cdot B(\mathbf{P})$ linear orders extending 
$\mathbf{P}$ and each $\mathbf{L}_{i,j}$, where
$i+k-1=j$. There are however $1 + 2 + \cdots + (m-\floor{(m+1)/2})$
ways of choosing a pair $i<j$ with $j-i \geq \floor{(m+1)/2}$.
Therefore the total number of linear orders 
extending $\mathbf{P}$ and each ${\mathbf{L}}_{i,j}\in{\mathcal{L}}_{i,j}$
for all $i<j$, will not exceed
\begin{eqnarray*}
\left[
\left(2+3+\cdots +\floor{\frac{m+1}{2}}\right) 
+ \left(1+2+\cdots +\left(m-\floor{\frac{m+1}{2}}\right)\right)\right]
\cdot B(\mathbf{P})\\
= \floor{\frac{(m-1)(m+3)}{4}}B(\mathbf{P}).
\end{eqnarray*}
Hence we have the theorem.
\end{proof}
{\sc Remark:} Considering the canonical interval order $\preceq$ 
on ${[m]\choose 2}$ in which $\{i_1,j_1\} \prec \{i_2,j_2\}$
iff $j_1 < i_2$, we see that all the $\floor{(m-1)(m+3)/4}$
2-sets $\{i,j\}\in {[m]\choose 2}$ 
with $i\in \{1,\ldots, \ceil{m/2}\}$ and $j\in \ceil{m/2},\ldots,m\}$ are incomparable. 
This means that the total number of linear orders in the proof of 
Theorem~\ref{thm:better-bounds}, that extend $\mathbf{P}$ and each 
${\mathbf{L}}_{i,j}\in{\mathcal{L}}_{i,j}$, cannot be reduced any further with 
the arguments presented there.

To better understand the asymptotic behavior of $\dim(\mathbf{P})$
of an $m$-partite poset $\mathbf{P}$, define $f(m)$ for each $m\geq 2$
by
\begin{equation}
\label{eqn:f(m)}
f(m) = \sup_{\mathbf{P}}\left\{\frac{\dim(\mathbf{P})}{B(\mathbf{P})}\right\},
\end{equation}
where the supremum is taken over all $m$-partite posets $\mathbf{P}$. 
By Theorem~\ref{thm:better-bounds} we therefore have that
$f(m)\leq \floor{(m-1)(m+3)/4}$.

For the lower bound of $f(m)$, we start with the following lemma.
\begin{lemma}
\label{lmm:complete-g}
Let $g,h,k\in \nats$ with $h,k\geq 2$ and $g\leq\min\{h,k\}$.
Let $M\subseteq [h]\times [k]$ be any matching of size $g$
between the columns and rows of $[h]\times [k]$.
For disjoint sets $X = \{x_1,\ldots,x_h\}$ and 
$Y = \{y_1,\ldots,y_k\}$ let $\mathbf{C}_{-g}(h,k)$ 
be the poset on $X\cup Y$ given by $x_i\prec y_j$ for all 
$(i,j)\in [h]\times [k]\setminus M$. Then 
$\dim(\mathbf{C}_{-g}(h,k)) = \max\{2,g\}$.
\end{lemma}
\begin{proof}
Assume $g\geq 2$. By a suitable permutation we may assume 
that $M = \{(1,1),\ldots,(g,g)\}$. Since the poset induced by 
$\{x_1,\ldots,x_g\}\cup \{y_1,\ldots,y_g\}$
is the standard example $\mathbf{S}_{2g}$ we have that 
$\dim(\mathbf{C}_{-g}(h,k)) \geq g$.

Let $L_x$ denote the linear order 
$x_1\prec x_3\prec x_4\prec\cdots\prec x_{h-1}\prec x_h\prec x_2$
and similarly let $L_y$ denote
$y_1\prec y_3\prec y_4\prec\cdots\prec y_{k-1}\prec y_k\prec y_2$.
If $i\in [h]$ then $L_x(\hat{i})$ denotes the linear order obtained
from $L_x$ by removing $x_i$ and similarly for $L_y(\hat{j})$.
For any linear order $L$ let $L^{\op}$ denote the opposite, or reverse,
linear order of $L$. In this case $\mathbf{C}_{-g}(h,k)$ is realized
by the following $g$ linear orders
\begin{eqnarray*}
L_x(\hat{1})^{\op}\prec y_1\prec x_1\prec L_y(\hat{1}), &  \\
L_x(\hat{\ell})\prec y_{\ell}\prec x_{\ell}\prec L_y(\hat{\ell})^{\op}
& \mbox{ for $\ell\in \{2,\ldots,g\}$. } 
\end{eqnarray*}
Hence $\dim(\mathbf{C}_{-g}(h,k)) \leq g$. The case $g=1$ gives in similar fashion 
$\dim(\mathbf{C}_{-1}(h,k)) = 2$.
\end{proof}
Note that $\mathbf{C}_{-g}(h,k))$ is the complete bipartite poset on $X$ and $Y$
except for the $g$ relations $x_i\prec y_j$ where $(i,j)\in M$.
\begin{theorem}
\label{thm:f(m)lower}
For $m\geq 2$ we have that $f(m)$ defined in (\ref{eqn:f(m)}) satisfies
\[
f(m)\geq \floor{\frac{m}{2}}\ceil{\frac{m}{2}} = \floor{\frac{m^2}{4}}.
\]
\end{theorem}
\begin{proof}
For $d,h,k\geq 2$ let $A = \{ x_{i,j} : (i,j)\in[dh]\times[dk]\}$
and $B = \{ y_{i,j} : (i,j)\in[dh]\times[dk]\}$ 
be two disjoint sets of $d^2hk$ elements each. 
Let $\mathbf{P} = (A\cup B;\preceq)$
be given by 
\[
x_{i_1,j_1}\prec y_{i_2,j_2} \Leftrightarrow (i_1,j_1) \neq (i_2,j_2).
\]
Here $\mathbf{P}$ is the standard example on $2d^2hk$ elements so
$\dim(\mathbf{P}) = d^2hk$. Let $X_1,\ldots,X_h,Y_1,\ldots,Y_k$ be given
by 
$X_p = \{x_{i,j} : (i,j)\in \{(p-1)d+1,\ldots, pd\}\times[dk]\}$
for each $p\in[h]$ and
$Y_q = \{y_{i,j} : (i,j)\in [dh]\times\{(q-1)d+1,\ldots, qd\}\}$
for each $q\in[k]$. This partition
of $A\cup B$ makes $\mathbf{P}$ 
into a $(h+k)$-partite poset $(X_1,\ldots,X_h,Y_1,\ldots,Y_k;\preceq)$.
We note that each of $X_p\cup X_q$ and $Y_p\cup Y_q$ is an antichain
in $\mathbf{P}$ of order dimension two. Since the sub-poset
of $\mathbf{P}$ induced by $X_p\cup Y_q$ is $\mathbf{C}_{-d^2}(d^2k,d^2h)$
we have by Lemma~\ref{lmm:complete-g} that
$B(\mathbf{P}) = d^2$. Hence we have 
\[
f(h+k) \geq \frac{\dim(\mathbf{P})}{B(\mathbf{P})} = \frac{d^2hk}{d^2} = hk.
\]
Putting $(h,k) = (n,n)$ on one hand and $(h,k) = (n,n+1)$ on the other
yields a lower bound for $f(m)$ both for even and odd $m$.
Hence, we have the theorem.
\end{proof}
Note that the example provided in the above proof of Theorem~\ref{thm:f(m)lower}
shows that both $\dim(\mathbf{P})$ and $B(\mathbf{P})$ can be arbitrarily
large.

By Theorems~\ref{thm:better-bounds} and~\ref{thm:f(m)lower}
we have the following.
\begin{corollary}
\label{cor:f(m)tight}
If $f(m)$ is the function from (\ref{eqn:f(m)}), then for all $m\geq 2$
we have
\[
\floor{\frac{m^2}{4}}\leq f(m)\leq \floor{\frac{(m-1)(m+3)}{4}}.
\]
\end{corollary}
By Corollary~\ref{cor:f(m)tight} we have
$\lim_{m\rightarrow\infty}f(m)/m^2 = 1/4$,
so the upper bound in Theorem~\ref{thm:better-bounds}
is asymptotically tight.

\section*{Acknowledgments}  
The author likes to thank Reinhard Laubenbacher for 
interesting discussions regarding applications of posets.
Also, sincere thanks to Walter D.~Morris for his helpful 
comments on the article.

\end{document}